\documentclass{birkjour}
\usepackage{cite}
\usepackage[colorlinks,linkcolor=blue,citecolor=green]{hyperref}
\usepackage{cleveref}
\usepackage{mathrsfs}
\usepackage{amsthm,amsmath,amssymb}

\crefname{equation}{}{}
\newtheorem{theorem}{Theorem}[section]
\newtheorem{lemma}[theorem]{Lemma}

\newtheorem{corollary}[theorem]{Corollary}

\theoremstyle{definition}
\newtheorem{definition}[theorem]{Definition}

\theoremstyle{remark}
\newtheorem{remark}[theorem]{Remark}

\numberwithin{equation}{section}

\begin{document}

%
%
%
%
%
%
%
%
%

\title[Liouville Theorem]
 {Boundary pointwise regularity and Liouville theorems for fully nonlinear equations on cones}

\author{Yuanyuan Lian}

\address{%
School of Mathematical Sciences\\
Shanghai Jiao Tong University\\
Shanghai, 200240\\
PR China}

\email{lianyuanyuan@sjtu.edu.cn; lianyuanyuan.hthk@gmail.com}

\thanks{This research is supported by the China Postdoctoral Science Foundation (Grant No.
2021M692086), the National Natural Science Foundation of China (Grant No. 12031012 and 11831003) and  the Institute of Modern Analysis-A Frontier Research Center of Shanghai.}


\date{}
\dedicatory{}

\begin{abstract}
In this paper, we prove a boundary pointwise regularity for fully nonlinear elliptic equations on cones. In addition, based on this regularity, we give simple proofs of the Liouville theorems on cones.

\end{abstract}

\maketitle
\section{Introduction}
\label{intro}
In this paper, we study the boundary pointwise regularity and Liouville theorems on cones for viscosity solutions of fully nonlinear uniformly elliptic equations
\begin{equation}\label{main}
F(D^2u)=0,
\end{equation}
where $F$ is a real fully nonlinear operator defined in $\mathcal{S}$ with ellipticity constants $\lambda$ and $\Lambda$. Here, $\mathcal{S}$ denotes the set of $n\times n$ symmetric matrices.
%

We are interested in understanding how solutions of \cref{main} that subject to Dirichlet boundary conditions behave at the vertex and the infinity of an unbounded cone. In other words, we investigate the boundary pointwise regularity and its connections to the Liouville theorems of \cref{main} on cones.

Pointwise regularity occupies an important position because it gives deeper understanding of the behaviour of the solution near the point concerned. Usually, the assumptions are weaker than that for local and global regularity and the pointwise regularity implies the local and global regularity. Many pointwise regularity results spring up since Caffarelli \cite{MR1005611} (see also \cite{MR1351007}) proved the interior pointwise $C^{1,\alpha}$ and $C^{2,\alpha}$ regularity for fully nonlinear elliptic equations. With respect to the boundary pointwise regularity, Silvestre and Sirakov \cite{MR3246039} proved the boundary pointwise $C^{1,\alpha}$ and $C^{2,\alpha}$ regularity on flat boundaries. Lian and Zhang \cite{MR4088470} obtained the boundary pointwise $C^{1,\alpha}$ and $C^{2,\alpha}$ regularity with respect to general smooth boundaries. However, the boundary pointwise regularity on domains with corners, e.g. cones, remains empty. In this paper, we will derive the boundary pointwise regularity on cones by a simple way.

The Liouville theorems have received many concerns because of its wide applications in several areas of mathematics such as differential geometry, geometric analysis, non-linear potential theory and PDEs. The Liouville theorems are fundamental to the classification of global solutions. The first Liouville theorem was presented by Liouville in 1844 and immediately proved by Cauchy\cite{FARINA200761}. Large amount research has been performed for different types of elliptic equations.

It is easy to show that any bounded harmonic function in $R^n$ must be a constant by the mean value property. For more complicated equations, the representation of solutions or energy estimates method are applied widely for the Liouville theorems.
For example, based on energy estimates, Gidas and Spruck \cite{MR615628} showed that any nonnegative solution $u$ is identically equal to $0$ for the semilinear equation $\Delta u+u^p=0$ ($1<p<n+2/n-2$). Another idea to obtain the Liouville theorem is that the rigidity of solutions can be derived from the symmetry of solutions. Chen and Li \cite{MR1121147} proved the same Liouville theorem by employing the moving plane method. For fully nonlinear elliptic equations, Armstrong, Sirakov and Smart \cite{MR2947535} proved the Liouville theorem on cones. Braga \cite{MR3735564} gave a simple proof of the Liouville theorem on half-spaces depending boundary Harnack inequality. In this paper, we presents a new and simple proof of the Liouville theorems on cones with the aid of the boundary pointwise regularity.

 We use standard notations in this paper. For $x\in R^n$, we may write $x=(x_1,...,x_n)=(x',x_n)$. Let $R^n_+=\left\{x\big|x_n>0\right\}$ denote the upper half-space. As usual, $|x|$ is the Euclidean norm of $x$ and the unit sphere is $S^{n-1}=\{x\big||x|=1\}$. Set $B_r(x_0)=\{x\big| |x-x_0|<r\}$ and $B_r=B_r(0)$. Additionally, the half ball is represented as $B_r^+(x_0)=B_r(x_0)\cap R^n_+$ and $B_r^+=B^+_r(0)$. Similarly, $T_r(x_0)\ =\{(x',0)\big| |x'-x_0'|<r\}$ and $T_r=T_r(0)$. The $\{e_i\}_{i=1}^{n}$ is the standard basis of $R^n$ and $\mathrm{dist}(A,B)$ denotes the distance between $A$ and $B$ for $A,B\subset R^n$.

\begin{definition}[Conical cone]\label{d-cone}
Let $\omega$ be a proper $C^2$ smooth subdomain of the
unit sphere $S^{n-1}$ ($n\geq 2$). Define the infinite conical cone
\begin{equation*}
C_{\omega}:= \left\{x\in R^n\big| |x|^{-1}x \in \omega \right\}.
\end{equation*}
In addition, we often consider the cone in a ball. Hence, we define $C^r_{\omega}=C_{\omega}\cap B_r$ and $(\partial C_{\omega})^r=\partial C_{\omega}\cap B_r$. Throughout this paper, we always assume that $e_n\in \omega$ and $0$ is the vertex of the cone without loss of generality.
\end{definition}

Armstrong, Sirakov and Smart \cite{MR2947535} proved the existence of nonnegative homogenous solutions on conical cones:

\begin{lemma}\label{l.psi-}
Assume that $F$ is positively $1$-homogeneous, i.e.,
\begin{equation}\label{struc-2}
F(tM)=tF(M),~~\forall~~t\geq 0, M\in \mathcal{S}.
\end{equation}
Then for any $\omega\subset S^{n-1}$, there exists a nonnegative solution $\Psi\in C(\bar C_{\omega})$ of
\begin{equation*}
  \left\{\begin{aligned}
    &F(D^2 u)=0&& ~~\mbox{in}~~C_{\omega};\\
    &u=0&& ~~\mbox{on}~~\partial C_{\omega}
  \end{aligned}\right.
\end{equation*}
with
\begin{equation*}
  \Psi(tx)=t^{\alpha_{\omega}}\Psi(x),~~~~\forall~~t>0,~~x\in C_{\omega},
\end{equation*}
where $\alpha_{\omega}>0$ depends only on $F$ and $\omega$.
\end{lemma}

\begin{remark}\label{re-Psi-}
Throughout this paper, we always assume (excepting the case of a half-space, see \Cref{sec:1}) that $F$ is positively $1$-homogeneous unless stated otherwise.
\end{remark}

Our main results are the following.

\begin{theorem}[\textbf{Pointwise $C^{\alpha_{\omega},\alpha}$ regularity}]\label{t.reg-F-cone}
Let $u$ be a viscosity solution of
\begin{equation}\label{e.reg-cone-1}
\left\{\begin{aligned}
&F(D^2u)=0&& ~~\mbox{in}~~C^1_{\omega};\\
&u=0&& ~~\mbox{on}~~(\partial C_{\omega})^1.
\end{aligned}\right.
\end{equation}
Then there exists a constant $a$ such that
\begin{equation}\label{e.reg-cone-res}
    |u(x)-a\Psi(x)|\leq C|x|^{\alpha}\Psi(x) \|u\|_{L^{\infty}(C^1_{\omega})} , ~~\forall ~x\in C^{\frac{1}{2}}_{\omega}
\end{equation}
and
\begin{equation*}
  |a|\leq C,
\end{equation*}
where $\alpha$ and $C$ depend only on $n,\lambda,\Lambda$ and $\omega$.
\end{theorem}

\begin{remark}\label{t.reg-cone-3}
\Cref{t.reg-F-cone} shows that the difference between $u$ and $a\Psi$ is controlled by an infinitesimal with order $\alpha$ higher than $\Psi$. Thus, we call $u\in C^{\alpha_{\omega},\alpha}(0)$. If $\Psi$ is a homogenous polynomial of degree $k$ (a nonnegative integer), we arrive at the classical boundary pointwise $C^{k,\alpha}$ regularity.
\end{remark}

\begin{remark}\label{t.reg-cone-1}
Interestingly, on the upper half-space $R^n_{+}$, $\Psi$ is exactly $x_n$ and $\alpha_{R^n_{+}}=1$ for every uniformly elliptic equation. Thus, \Cref{t.reg-F-cone} reduces to the well-kown boundary $C^{1,\alpha}$ regularity for equations on half-spaces. This was first proved by Krylov \cite{MR688919} and further simplified by Caffarelli \cite[Theorem 9.31]{MR1814364}.
\end{remark}

\begin{remark}\label{t.reg-cone-2}
By checking the definition of $\alpha_{\omega}$ (see \cite[3.2]{MR2947535}), $\alpha_{\omega'}>\alpha_{\omega}$ if $\omega'\subset \omega$. That is, the regularity at the vertex of the cone is higher if $\omega$ is smaller.
\end{remark}

\begin{remark}\label{re1.1}
The boundary pointwise regularity also holds for general domains with corners, such as Lipschitz domains. The general cases:
\begin{equation*}
\left\{\begin{aligned}
&F(D^2u,Du,u,x)=f&& ~~\mbox{in}~~\Omega;\\
&u=g&& ~~\mbox{on}~~\partial \Omega
\end{aligned}\right.
\end{equation*}
will be treated in a future work.
\end{remark}

If $F$ is the Laplace operator, $\alpha_{\omega}$ can be expressed explicitly as
\begin{equation}\label{e1.1}
\alpha_{\omega}=\frac{1}{2}\left(\sqrt{(n-2)^2+4\lambda_1(\omega)}-(n-2)\right),
\end{equation}
where $\lambda_1(\omega)$ the first eigenvalue of the Dirichlet problem for the spherical Laplacian in $\omega$ (see \cite{MR2904138} for example). The $\Psi$ is exactly the eigenfunction corresponding to $\lambda_1(\omega)$. Hence, we have the following interesting corollary.

\begin{corollary}\label{t.laplace}
Let $u$ be a viscosity solution of
\begin{equation}\label{e.reg-lap-cone-1}
\left\{\begin{aligned}
&\Delta u=0&& ~~\mbox{in}~~C^1_{\omega};\\
&u=0&& ~~\mbox{on}~~\partial C^1_{\omega}.
\end{aligned}\right.
\end{equation}
Then $u\in C^{\alpha_{\omega},\alpha}(0)$, i.e., there exists a constant $a$ such that
\begin{equation}\label{e.reg-lap-cone-res}
    |u(x)-a\Psi(x)|\leq C|x|^{\alpha}\Psi(x) \|u\|_{L^{\infty}(C^1_{\omega})} , ~~\forall ~x\in C^{\frac{1}{2}}_{\omega}
\end{equation}
and
\begin{equation*}
  |a|\leq C,
\end{equation*}
where $\alpha$ and $C$ depend only on $n,\lambda,\Lambda$ and $\omega$ and $\alpha_{\omega}$ is defined as in \cref{e1.1}.

In particular, if $n=2$,
\begin{equation}\label{e1.2}
\alpha_{\omega}=\sqrt{\lambda_1(\omega)}=\frac{\pi}{\theta},
\end{equation}
where $\theta$ is the central angle of the cone $C_{\omega}$ (reduced to a circular sector). Hence, $u\in C^{\pi/\theta,\alpha}(0)$ for some $0<\alpha<1$ depending only on $\theta$.
\end{corollary}

\begin{remark}\label{rk.t.laplace-1}
 In the interesting case of $n=2$, we have $u\in C^{1,\alpha}(0)$ on the half-space and $u\in C^{2,\alpha}(0)$ on the first quadrant etc. In addition, $u\in C^{\alpha}(0)$ with $\alpha>1/2$ for any cone since the central angle $\theta<2\pi$.
\end{remark}

Based on the boundary pointwise regularity, we can prove the following Liouville theorems.
\begin{theorem}[Liouville theorem-type I]\label{t.reg-co}
Let $u$ be a viscosity solution of
\begin{equation*}\label{e.ful-reg-co}
\left\{\begin{aligned}
&F(D^2u)=0&& ~~\mbox{in}~~C_{\omega};\\
&u=0&& ~~\mbox{on}~~\partial C_{\omega}.
\end{aligned}\right.
\end{equation*}
If $u(x)=O(\Psi(x))$ as $x\rightarrow \infty$,
\begin{equation*}
  u\equiv \frac{u(e_n)}{\Psi(e_n)}\Psi~~\mbox{ in}~~C_{\omega},
\end{equation*}
In particular, if $u(x)=o(\Psi(x))$ as $x\rightarrow \infty$,
\begin{equation*}
  u\equiv 0~~\mbox{ in}~~C_{\omega}.
\end{equation*}
\end{theorem}

\begin{theorem}[\textbf{Liouville theorem-type II}]\label{t.ful-cone}
Let $u\geq 0$ be a viscosity solution of
\begin{equation*}\label{e.ful-cone}
\left\{\begin{aligned}
&F(D^2u)=0&& ~~\mbox{in}~~C_{\omega};\\
&u=0&& ~~\mbox{on}~~\partial C_{\omega}.
\end{aligned}\right.
\end{equation*}
Then
\begin{equation}\label{e.ful-cone-Psi}
  u\equiv\frac{u(e_n)}{\Psi(e_n)}\Psi ~~~\mbox{ in}~~ C_{\omega}.
\end{equation}
\end{theorem}



This paper is organized as follows. We first consider the special cone, i.e., $R^n_+$ in \Cref{sec:1}. In \Cref{sec:3}, we prove the pointwise regularity at the vertex of a conical cone and then obtain the Liouville theorems on cones in the similar way as on half-spaces. Throughout this paper, we use the letter $C$ denote a positive constant depending only on the dimension $n$, the ellipticity constants $\lambda$ and $\Lambda$ and $\omega$ and we call $C$ a universal constant.

\section{Liouville theorems on half-spaces}\label{sec:1}
To show the main idea in a clear manner, we prove the Liouville theorems on $R^n_+$ based on the boundary pointwise regularity in this section. In addition, the homogenous condition \cref{struc-2} is not needed and hence we don't assume \cref{struc-2} throughout this section. Another reason for considering $R^n_+$ and general cones separately is that $\Psi\equiv Cx_n$ on $R^n_+$ for any elliptic operator $F$. Indeed, we prove the Liouville theorems for the Pucci class (see the main results below). Furthermore, we can obtain a higher order Liouville theorem on $R^n_+$ based on the same idea (see \Cref{t.ful-2}).



We first introduce two lemmas (see \cite[Lemma 2.11 and Lemma 2.12]{WLZ}).
\begin{lemma}[\textbf{Hopf lemma}]\label{l.bas}
Let $u(e_n/2)=1$ and $u\geq 0$ be a viscosity solution of
\begin{equation}
  \left\{\begin{aligned}
    &u\in S(\lambda,\Lambda,0)&& ~~\mbox{in}~~B_1^+;\\
    &u=0&& ~~\mbox{on}~~T_1.
  \end{aligned}\right.
\end{equation}
Then
\begin{equation}\label{e.bas}
u(x)\geq Cx_n ~~\mbox{in}~~B^+_{1/2},
\end{equation}
where $C$ is universal.
\end{lemma}

\begin{remark}\label{re2.1}
The $S(\lambda,\Lambda,0)$ denotes the Pucci class. For more details about its properties, we refer to \cite{MR1351007}.
\end{remark}
%

\begin{theorem}[\textbf{Boundary pointwise $C^{1,\alpha}$ regularity}]\label{t.reg-S}
Let $u$ be a viscosity solution of
\begin{equation}\label{e.reg-1}
\left\{\begin{aligned}
&u\in S(\lambda,\Lambda,0)&& ~~\mbox{in}~~B_1^+;\\
&u=0&& ~~\mbox{on}~~T_1.
\end{aligned}\right.
\end{equation}
Then $u\in C^{1,\alpha}(0)$, i.e., there exists a constant $a$ such that
\begin{equation}\label{e.reg-res}
    |u(x)-ax_n|\leq C|x|^{\alpha}x_n \|u\|_{L^{\infty}(B_1^+)}, ~~\forall ~x\in B_{1/2}^+
\end{equation}
and
\begin{equation*}
  |a|\leq C\|u\|_{L^{\infty}(B_1^+)},
\end{equation*}
where $0<\alpha<1$ and $C$ are universal.
\end{theorem}

From the boundary pointwise regularity, we can obtain the Liouville theorem-type I immediately:
\begin{theorem}[Liouville theorem-type I]\label{co2.1}
Let $u$ be a viscosity solution of
\begin{equation}\label{e.reg-1}
\left\{\begin{aligned}
&u\in S(\lambda,\Lambda,0)&& ~~\mbox{in}~~R_+^n;\\
&u=0&& ~~\mbox{on}~~\partial R_+^n.
\end{aligned}\right.
\end{equation}
If $u=O(x_n)$ as $x\rightarrow \infty$,
\begin{equation*}
  u\equiv u(e_n)x_n~~\mbox{ in}~~R_+^n,
\end{equation*}
In particular, if $u=o(x_n)$ as $x\rightarrow \infty$,
\begin{equation*}
  u\equiv 0~~\mbox{ in}~~R_+^n.
\end{equation*}
\end{theorem}
\proof A scaling version of \Cref{t.reg-S} gives that for any $R>0$,
\begin{equation}\label{e.t.ful-reg-0}
|u(x)-ax_{n}|\leq C\frac{|x|^{1+\alpha}}{R^{1+\alpha}}\cdot \|u\|_{L^{\infty}(B_R^+)},
~~\forall~x\in B^+_{R/2},
\end{equation}
where
\begin{equation*}
  |a|\leq C\|u\|_{L^{\infty}(B_R^+)}/R.
\end{equation*}

Note that $u=O(x_n)$. Fix $x\in R^n_+$ and let $R\rightarrow \infty$ in \cref{e.t.ful-reg-0}. Then $u(x)= ax_n$  and hence
\begin{equation*}
  u\equiv ax_n~\mbox{ in }~R^n_+.
\end{equation*}
Obviously, $a=u(e_n)$.

If $u=o(x_n)$, $a$ must be zero. That is, $u\equiv 0$ in $R^n_+$.~\qed~\\

To prove the Liouville theorem-type II, we need the following Carleson type estimate.
\begin{lemma}[\textbf{Carleson type estimate}]\label{l.cal}
Let $u(e_n/2)=1$ and $u\geq 0$ be a viscosity solution of
\begin{equation*}\label{e.cal-1}
\left\{\begin{aligned}
&u\in S(\lambda,\Lambda,0)&& ~~\mbox{in}~~B_1^+;\\
&u=0&& ~~\mbox{on}~~T_1.
\end{aligned}\right.
\end{equation*}
Then
\begin{equation}\label{e.cal-2}
  \|u\|_{L^{\infty}(B^+_{1/2})}\leq C,
\end{equation}
where $C$ is universal.
\end{lemma}

\proof Let $v$ be the zero extension of $u$ to the whole $B_1$. Then $v$ satisfies (see \cite[Proposition 2.8]{MR1351007})
\begin{equation*}
 v\in \underline{S}(\lambda,\Lambda,0)~~\mbox{ in}~~B_1.
\end{equation*}
By the local maximum principle (see \cite[Theorem 4.8]{MR1351007}), for any $p>0$,
\begin{equation*}\label{cal-lp}
  \|u\|_{L^{\infty}(B_{1/2}^+)}\leq \|v\|_{L^{\infty}(B_{1/2})}
  \leq C_0 \|v\|_{L^{p}(B_{3/4})}= C_0\|u\|_{L^{p}(B^+_{3/4})},
\end{equation*}
where $C_0$ depends only on $n,\lambda,\Lambda$ and $p$. Thus, we only need to prove that for some universal constant $p>0$,
\begin{equation*}\label{e2.1}
\|u\|_{L^{p}(B^+_{3/4})}\leq C.
\end{equation*}

Fix $x_0 \in \bar{B}^+_{3/4}$ and denote $L(r,R)=\left\{(x',x_n): x'=x'_0,r\leq x_n\leq R\right\}$ for $0<r<R$. From $u(e_n/2)=1$ and the Harnack inequality (see \cite[Theorem 4.3]{MR1351007}),
\begin{equation*}
  \sup_{L(1/2,3/4)} u\leq C.
\end{equation*}
By the Harnack inequality again,
\begin{equation*}
  \sup_{L(1/4,1/2)} u\leq C_1
\end{equation*}
where $C_1$ is universal. Similarly,
\begin{equation*}
  \sup_{L(1/2^k,1/2^{k-1})} u\leq C_1^{k-1},~~~~\forall~~k\geq 2.
\end{equation*}

Take $k\geq 2$ such that $1/2^k\leq x_{0,n}\leq 1/2^{k-1}$. Then
\begin{equation}\label{e2.5}
  u(x_0)\leq \sup_{L(1/2^k,1/2^{k-1})} u \leq C_1^{k-1}=C(2^{-k})^{-q}\leq C x_{0,n}^{-q},
\end{equation}
where $q>0$ is universal. Hence,
\begin{equation*}
\|u\|_{L^p(B^+_{3/4})}\leq C
\end{equation*}
for some universal $p>0$. \qed~\\
\begin{remark}\label{re2.2}
To obtain the Carleson type estimate by combining the interior Harnack inequality and zero extension is inspired by the work of De Silva and Savin \cite{MR4093736}.
\end{remark}

Now, we can prove the Liouville theorem-type II on the half-space.
\begin{theorem}[\textbf{Liouville theorem-type II}]\label{t.ful-1}
Let $u\geq 0$ be a viscosity solution of
\begin{equation}\label{e.ful-1}
\left\{\begin{aligned}
&u\in S(\lambda,\Lambda,0)&& ~~\mbox{in}~~R^n_+;\\
&u=0&& ~~\mbox{on}~~\partial R^n_+.
\end{aligned}\right.
\end{equation}
Then
\begin{equation*}\label{e.ful-1-xn}
  u(x)=u(e_n)x_n, ~~\forall ~x\in R^n_+.
\end{equation*}
\end{theorem}

\proof By \Cref{t.reg-S} (in fact the scaling version of \Cref{t.reg-S}), there exists a constant $a$ such that for any $R>0$,
\begin{equation}\label{e.t.ful-reg}
|u(x)-ax_{n}|\leq C \frac{|x|^{1+\alpha}}{R^{1+\alpha}}\cdot \|u\|_{L^{\infty}(B_R^+)},
~~\forall~x\in B^+_{R/2}.
\end{equation}
The Carleson estimate (\Cref{l.cal}) implies
\begin{equation}\label{e2.2}
  \|u\|_{L^{\infty}(B_R^+)}\leq Cu(Re_n/2).
\end{equation}
From the Hopf lemma (\Cref{l.bas}),
\begin{equation*}
  u(x)\geq C u(Re_n/2)\cdot \frac{x_n}{R},~~\forall~~x\in B^+_{R/2}.
\end{equation*}
By setting $x=e_n/2$ and we have
\begin{equation}\label{e2.3}
  u(Re_n/2)\leq C u(e_n/2) R.
\end{equation}

By combining \cref{e.t.ful-reg}, \cref{e2.2} and \cref{e2.3}, we have
\begin{equation*}
|u(x)-ax_{n}|\leq Cu(e_n/2)\cdot \frac{|x|^{1+\alpha}}{R^{\alpha}},
~~\forall~x\in B^+_{R/2}.
\end{equation*}
Fix any $x\in R^n_+$ and let $R\rightarrow \infty$. Then $u(x)=ax_n$ and $a=u(e_n)$. Hence,
\begin{equation*}
  u\equiv u(e_n)x_n~~~~\mbox{  in }~R^n_+.
\end{equation*}
\qed~\\

\begin{remark}\label{r-1.1}
If we consider the Liouville theorem in the whole space , $u\equiv 0$ can be obtained by the Harnack inequality. However, on a half-space, $u$ has a linear growth away from the flat boundary according to the Hopf lemma. Hence, $u$ is a linear function in $R^n_+$.
\end{remark}

\begin{remark}\label{r-1.2}
The method of proving the Liouville theorem based on the boundary pointwise regularity is shown clearly in above proof. The boundary pointwise regularity gives the estimate of the error between the solution and a linear function (see \cref{e.t.ful-reg}). If $\|u\|_{L^{\infty}(B_R^+)}$ grows linearly as $R\rightarrow \infty$, we arrive at the Liouville theorem. The Hopf lemma guarantees that $u(Re_n/2)$ grows linearly. Finally, the Carleson type estimate provides a bridge between $\|u\|_{L^{\infty}(B_R^+)}$ and $u(Re_n/2)$.
\end{remark}

\begin{remark}\label{re2.3}
\Cref{t.ful-1} has been proved by Armstrong, Sirakov and Smart \cite{MR2947535} and Braga gave a simplified proof. Our proof is simpler and clearer by comparing with that of \cite{MR3735564}.
\end{remark}

Since $R^n_+$ is a special cone, we can prove a higher order Liouville theorem with aid of the boundary pointwise regularity. The next lemma concerns the boundary pointwise $C^{2,\alpha}$ regularity (see Silvestre and Sirakov and \cite{MR4088470}).

\begin{lemma}\label{l.reg-F}
Let $u$ be a viscosity solution of
\begin{equation*}
\left\{\begin{aligned}
&F(D^2u)=1&& ~~\mbox{in}~~B_1^+;\\
&u=0&& ~~\mbox{on}~~T_1.
\end{aligned}\right.
\end{equation*}
Then $u\in C^{2,\alpha}(0)$, i.e., there exists a polynomial $P$ of degree $2$ such that
\begin{equation*}\label{e.reg-1}
  |u(x)-P(x)|\leq C|x|^{2+\alpha}(\|u\|_{L^{\infty }(B_1^+)}+1+|F(0)|), ~~\forall ~x\in B_{1/2}^+
\end{equation*}
and
\begin{equation*}
  F(D^2P)=1,
\end{equation*}
where $0<\alpha<1$ and $C$ are universal constants. Moreover, $P$ can be written as
\begin{equation}\label{e2.4}
  P(x)=\sum_{i=1}^{n}a_{in}x_ix_n+bx_n,
\end{equation}
where $a_{in}$ and $b$ are constants.
\end{lemma}

Next, we prove a higher order Liouville theorem on $R^n_+$.
\begin{theorem}\label{t.ful-2}
Let $u\geq 0$ be a viscosity solution of
\begin{equation}\label{e.ful-2}
\left\{\begin{aligned}
&F(D^2u)=1&& ~~\mbox{in}~~R^n_+;\\
&u=0&& ~~\mbox{on}~~\partial R^n_+;\\
&D u=0&& ~~\mbox{on}~~\partial R^n_+.
\end{aligned}\right.
\end{equation}
Then
\begin{equation}\label{e.ful-2-xn2}
  u(x)=ax_n^2, ~~\forall ~x\in R^n_+,
\end{equation}
where $a$ depends only on $F$.
\end{theorem}

\proof  Let $v$ be the zero extension of $u$ to the whole space $R^n$. Since $u\in C_{loc}^{2,\alpha}(\bar R^n_+)$ and $Du=0$ on $\partial R^n_+$, $v\in C_{loc}^{1,1}(R^n)$. Hence, $v$ is a nonnegative strong solution of
\begin{equation}\label{e.t.ful2-ext}
\begin{aligned}
 &F(D^2v)=\chi_{R^n_+}&&~~\mbox{in}~~R^n.
  \end{aligned}
\end{equation}

By the Harnack inequality and $v(0)=0$, for any $R>0$,
\begin{equation}\label{e.t.ful2-1}
\|u\|_{L^{\infty}(B_R^+)}=\sup_{B_R} v
\leq C(\inf_{B_R} v+R^2\|\chi_{R^n_+}\|_{L^{\infty}(B_R)}+R^2|F(0)|)\leq CR^2.
\end{equation}
For any fixed $x \in R^n_+$, by taking $R>2|x|$ and \Cref{l.reg-F}, there exists a polynomial $P$ in the form \cref{e2.4} such that
\begin{equation}\label{e.t.ful2-2}
\begin{aligned}
|u(x)-P(x)|
&\leq C \frac{|x|^{2+\alpha}}{R^{2+\alpha}}(\|u\|_{L^{\infty}(B_R^+)}+R^2+R^2|F(0)|)\\
&=C \frac{|x|^{2+\alpha}}{R^{2+\alpha}}(\|u\|_{L^{\infty}(B_R^+)}+R^2)
\end{aligned}
\end{equation}

By combining \cref{e.t.ful2-1} and \cref{e.t.ful2-2} and letting $R\rightarrow \infty$, we have
\begin{equation*}
  u\equiv P~\mbox{ in }~R^n_+.
\end{equation*}
Since $Du=0$ on $\partial R^n_+$, \cref{e.ful-2-xn2} holds. Clearly, $a$ is uniquely determined by $F(D^2P)=1$.~\qed~\\



\section{Boundary pointwise regularity and Liouville theorems on cones}\label{sec:3}

In this section, we first give the proof of the boundary pointwise regularity \Cref{t.reg-F-cone} and then prove the Liouville theorems \Cref{t.reg-co} and \Cref{t.ful-cone}. Similar to the Liouville theorems on half-spaces, we can obtain a control of the error between the solution and $\Psi$ by the boundary pointwise regularity at the vertex of the cone. Next, we use the Hopf lemma and the Carleson estimate to get a control of $u$. Then, the Liouville type theorem on cones can be derived.

The following Hopf lemma and boundary Lipschitz regularity can be proved easily by constructing proper barriers based on the interior and the exterior sphere conditions.

\begin{lemma}[Hopf lemma]\label{l.bas.C2}
Let $u(e_n/2)=1$ and $u\geq 0$ be a viscosity solution of
\begin{equation}
  \left\{\begin{aligned}
    &u\in S(\lambda,\Lambda,0)&& ~~\mbox{in}~~\Omega \cap B_1;\\
    &u=0&& ~~\mbox{on}~~\partial \Omega\cap B_1.
  \end{aligned}\right.
\end{equation}
Suppose that $0\in \partial \Omega$ and $\Omega$ satisfies the interior sphere condition at $0$, i.e.,
$B(r_0e_n,r_0)\subset \Omega$ for some $0<r_0<1/2$.

Then
\begin{equation}\label{e.bas.C2}
u(x)\geq Cx_n ,~~\forall ~x\in\left\{(0,x_n):0<x_n<r_0/2\right\},
\end{equation}
where $C>0$ depends only on $n,\lambda,\Lambda$ and $r_0$.
\end{lemma}

\begin{lemma}[Boundary Lipschitz regularity]\label{l.reg-C2}
Let $u$ be a viscosity solution of
\begin{equation*}
\left\{\begin{aligned}
&u\in S(\lambda,\Lambda,0)&& ~~\mbox{in}~~\Omega\cap B_1;\\
&u=0&& ~~\mbox{on}~~\partial \Omega\cap B_1.
\end{aligned}\right.
\end{equation*}
Suppose that $0\in \partial \Omega$ and $\Omega$ satisfies the exterior sphere condition at $0$, i.e.,
$B(-r_0e_n,r_0)\subset \Omega^c$ for some $0<r_0<1/2$.

Then $u$ is $C^{0,1}$ at $0$ and there exists a constant $a$ such that
\begin{equation}\label{e.l.C2-1}
  |u(x)|\leq C |x||u\|_{L^{\infty }(\Omega\cap B_1)}, ~\forall ~x\in \Omega\cap B_{r_0/2},
\end{equation}
where $C$ depends only on $n, \lambda, \Lambda$ and $r_0$.
\end{lemma}

\begin{remark}\label{rk.C2}
The geometrical conditions in \Cref{l.bas.C2} and \Cref{l.reg-C2} can be relaxed. In fact, \Cref{l.bas.C2} holds if $\Omega$ satisfies the interior $C^{1,\mathrm{Dini}}$ condition at $0$ and \Cref{l.reg-C2} holds if $\Omega$ satisfies the exterior $C^{1,\mathrm{Dini}}$ condition at $0$ (see \cite{lian2018boundary}).
\end{remark}



%
%

Next, we prove a Hopf lemma on cones, which is a simple consequence of the comparison principle.
\begin{lemma}[\textbf{Hopf lemma on cones}]\label{l.bas-cone}
Let $u(e_n/2)=1$ and $u\geq 0$ satisfy
\begin{equation*}
  \left\{\begin{aligned}
    &F(D^2u)=0&& ~~\mbox{in}~~C^1_{\omega};\\
    &u=0&& ~~\mbox{on}~~(\partial C_{\omega})^1.
  \end{aligned}\right.
\end{equation*}
Then
\begin{equation}\label{e.bas}
u\geq C\Psi ~~\mbox{  in}~~C^{1/2}_{\omega},
\end{equation}
where $C$ is universal.
\end{lemma}

\proof Since $\omega\in C^2$, $\partial C_{\omega}\in C^2$ excepting the origin. Hence,
$C_{\omega}\cap \left(B_{3/4}\backslash B_{1/4}\right)$ satisfies the uniform interior and exterior sphere condition with some radius $0<r_0<1/4$ (depending only on $\omega$) at any $x\in \partial C_{\omega}\cap \left(B_{3/4}\backslash B_{1/4}\right)$.

For any $x_0\in \partial B_{1/2}\cap C_{\omega}$, if $\mathrm{dist}(x_0,\partial C_{\omega})\geq r_0/2$, by the Harnack inequality
\begin{equation}\label{e3.1}
u(x_0)\geq Cu(e_n/2)=C\geq C\Psi(e_n/2)\geq C\Psi(x_0).
\end{equation}
If $\mathrm{dist}(x_0,\partial C_{\omega})< r_0/2$,
take $y_0\in \partial C_{\omega}$ with $|x_0-y_0|=d(x_0,\partial C_{\omega})$.
By \Cref{l.bas.C2} (up to a proper changing of coordinate system),
\begin{equation*}
  u(x_0)\geq C|x_0-y_0|.
\end{equation*}
In addition, from \Cref{l.reg-C2},
\begin{equation*}
  \Psi(x_0)\leq C|x_0-y_0|.
\end{equation*}
Hence,
\begin{equation}\label{e3.2}
u(x_0)\geq C\Psi(x_0).
\end{equation}

By combining \cref{e3.1} and \cref{e3.2}, we have
\begin{equation*}
  u\geq C\Psi ~~\mbox{on}~~\partial B_{1/2}\cap C_{\omega}.
\end{equation*}
Then according to the comparison principle (note that $u=\Psi=0$ on $(\partial C_{\omega})^{1/2}$),
\begin{equation*}
u\geq C\Psi~~\mbox{ in}~~C_{\omega}^{1/2}.
\end{equation*}~\qed\\

\begin{lemma}[\textbf{Boundary Lipschitz regularity on cones}]\label{l.blip-cone}
Let $u$ be a viscosity solution of
\begin{equation*}
  \left\{\begin{aligned}
    &F(D^2u)=0&& ~~\mbox{in}~~C^1_{\omega};\\
    &u=0&& ~~\mbox{on}~~(\partial C_{\omega})^1.
  \end{aligned}\right.
\end{equation*}
Then
\begin{equation}\label{e.blip}
|u(x)|\leq C\Psi(x) \|u\|_{L^{\infty}(C^1_{\omega})},~~~\forall ~x\in C^{1/2}_{\omega},
\end{equation}
where $C$ is universal.
\end{lemma}

\proof Without loss of generality, we assume that $\|u\|_{L^{\infty}(C^1_{\omega})}=1$. Similar to the proof of \Cref{l.bas-cone}, $C_{\omega}\cap \left(B_{3/4}\backslash B_{1/4}\right)$ satisfies the uniform interior and exterior sphere condition with some radius $0<r_0<1/4$ (depending only on $\omega$) at any $x\in \partial C_{\omega}\cap \left(B_{3/4}\backslash B_{1/4}\right)$.

For any $x_0\in \partial B_{1/2}\cap C_{\omega}$, if $\mathrm{dist}(x_0,\partial C_{\omega})\geq r_0/2$, by the Harnack inequality and  $\|u\|_{L^{\infty}(C^1_{\omega})}=1$,
\begin{equation}\label{e3.3}
\Psi(x_0)\geq\Psi(e_n/2)\geq C\geq Cu(x_0).
\end{equation}
If $\mathrm{dist}(x_0,\partial C_{\omega})< r_0/2$,
take $y_0\in \partial C_{\omega}$ with $|x_0-y_0|=d(x_0,\partial C_{\omega})$.
By \Cref{l.bas.C2} (up to a proper changing of coordinate system),
\begin{equation*}
  \Psi(x_0)\geq C|x_0-y_0|.
\end{equation*}
In addition, from \Cref{l.reg-C2},
\begin{equation*}
  u(x_0)\leq C|x_0-y_0|.
\end{equation*}
Hence,
\begin{equation}\label{e3.4}
\Psi(x_0)\geq Cu(x_0).
\end{equation}

By combining \cref{e3.3} and \cref{e3.4}, we have
\begin{equation*}
  \Psi\geq Cu ~~\mbox{ on}~~\partial B_{1/2}\cap C_{\omega}.
\end{equation*}
Then from the comparison principle,
\begin{equation*}
\Psi\geq Cu~~\mbox{ in}~~C_{\omega}^{1/2}.
\end{equation*}
Finally, consider $-u$ in above proof and then
\begin{equation*}
\Psi\geq -Cu~~\mbox{ in}~~C_{\omega}^{1/2}.
\end{equation*}
Hence, \cref{e.blip} holds.~\qed\\

We also need the following lemma, which is motivated by \cite[Lemma 2.54]{MR2947535}. For $\omega\subset S^{n-1}$ and $0<r<R$, denote
$E(\omega,r,R)=C_{\omega}\cap \left(B_R\backslash \bar B_r\right)$.
\begin{lemma}\label{le.3.1}
Let $\omega'\subset\subset \omega$ and $u$ satisfy
\begin{equation*}
  \left\{\begin{aligned}
    &M^-(D^2u)\leq 0&& ~~\mbox{in}~~E(\omega,1/8,1);\\
    &u\geq 1&& ~~\mbox{in}~~ E(\omega',1/4,3/4);\\
    &u\geq 0&& ~~\mbox{on}~~\partial C_{\omega}\cap (B_1\backslash \bar B_{1/8});\\
    &u\geq -\varepsilon_0&& ~~\mbox{on}~~\partial (B_1\backslash \bar B_{1/8})\cap C_{\omega},
  \end{aligned}\right.
\end{equation*}
where $0<\varepsilon_0<1$ depending only on $n,\lambda,\Lambda, \omega$ and $\omega'$ is small enough. Then
\begin{equation*}
  u\geq 0~~\mbox{ in}~~E(\omega,1/4,3/4).
\end{equation*}
\end{lemma}
\proof Let $v$ and $w$ be viscosity solutions of
\begin{equation*}
  \left\{\begin{aligned}
    &M^-(D^2v)= 0&& ~~\mbox{in}~~E(\omega,1/8,1)\backslash E(\omega',1/4,3/4);\\
    &v= 1&& ~~\mbox{in}~~ \partial E(\omega',1/4,3/4);\\
    &v= 0&& ~~\mbox{on}~~\partial E(\omega,1/8,1)
  \end{aligned}\right.
\end{equation*}
and
\begin{equation*}
  \left\{\begin{aligned}
    &M^+(D^2w)= 0&& ~~\mbox{in}~~E(\omega,1/8,1)\backslash E(\omega',1/4,3/4);\\
    &w= 0&& ~~\mbox{in}~~ \partial E(\omega',1/4,3/4);\\
    &w= 1&& ~~\mbox{on}~~\partial E(\omega,1/8,1)
  \end{aligned}\right.
\end{equation*}
respectively.

As before, $C_{\omega}$ satisfies the uniform interior and exterior sphere conditions with some radius $0<r_0<1/4$ (depending only on $\omega$) at any $x\in \partial C_{\omega}\cap \left(B_{3/4}\backslash B_{1/4}\right)$. For any $x_0\in E(\omega,1/8,1)\backslash E(\omega',1/4,3/4)$, if $\mathrm{dist}(x_0, \partial C_{\omega})\geq r_0/2$, by the strong maximum principle and noting $w(x_0)\leq 1$,
\begin{equation}\label{e3.6}
v(x_0)\geq C\geq Cw(x_0),
\end{equation}
where $C$ depends only on $n,\lambda,\Lambda,\omega$ and $\omega'$. If $\mathrm{dist}(x_0, \partial C_{\omega})< r_0/2$, by the Hopf lemma and the boundary Lipschitz regularity,
\begin{equation}\label{e3.7}
v(x_0)\geq C\mathrm{dist}(x_0, \partial C_{\omega})\geq Cw(x_0).
\end{equation}
Hence,
\begin{equation*}
v\geq Cw ~~\mbox{in}~~E(\omega,1/8,1)\backslash E(\omega',1/4,3/4).
\end{equation*}

Note that if $\varepsilon_0$ is small enough,
\begin{equation*}
  \left\{\begin{aligned}
    &M^-(D^2(v-Cw))\geq 0&& ~~\mbox{in}~~E(\omega,1/8,1)\backslash E(\omega',1/4,3/4);\\
    &v-Cw\leq u && ~~\mbox{in}~~ \partial E(\omega',1/4,3/4);\\
    &v-Cw\leq u&& ~~\mbox{on}~~\partial E(\omega,1/8,1).
  \end{aligned}\right.
\end{equation*}
Thus,
\begin{equation*}
  u\geq v-Cw\geq 0~~\mbox{in}~~E(\omega,1/8,1)\backslash E(\omega',1/4,3/4)
\end{equation*}
and therefore the conclusion holds.~\qed\\

Now, we can prove the ``$C^{\alpha_{\omega},\alpha}$ regularity'' at the vertex of a conical cone.

\noindent \textbf{Proof of \Cref{t.reg-F-cone}.} Without loss of generality, we assume that $\|u\|_{L^{\infty}(C^1_{\omega})}=1$. To obtain\cref{e.reg-cone-res}, we only need to prove the following: there exist a non-increasing sequence $\{a_k\}$ ($k\geq 0$) and a non-decreasing sequence $\{b_k\}$ ($k\geq 0$) such that
\begin{equation}\label{e3.5}
 |a_0|\leq C,~~ |b_0|\leq C
\end{equation}
and for all $k\geq 1$,
\begin{equation}\label{e.reg-cone-M}
\begin{aligned}
  &b_k\Psi
   \leq u\leq a_k\Psi~~~~\mbox{ in}~~~~C^{2^{-k}}_{\omega},\\
  &0\leq a_k-b_k\leq (1-\mu)(a_{k-1}-b_{k-1}),
\end{aligned}
\end{equation}
where $C$ and $0<\mu<1$ are universal.

We prove the above by induction. From \Cref{l.blip-cone},
\begin{equation*}
  -C\Psi\leq u\leq C\Psi~~\mbox{ in}~~C^{1/2}_{\omega}.
\end{equation*}
Thus, for $k=1$, by taking $a_0=2C, b_0=-2C$ and $a_1 = C, b_1 = -C$, \cref{e.reg-cone-M} holds
where $0<\mu<1/2$ is to be specified later. Assume \cref{e.reg-cone-M} holds for $k$, we need to prove that it holds for $k+1$.

Since \cref{e.reg-cone-M} holds for $k$, there are two possible cases:
\begin{equation*}
\begin{aligned}
&\mbox{\textbf{Case 1:}}~~&&~~u(2^{-k-1}e_n)\geq \frac{a_k+b_k}{2}\cdot \Psi(2^{-k-1}e_n);\\
&\mbox{\textbf{Case 2:}}~~&&~~u(2^{-k-1}e_n)< \frac{a_k+b_k}{2}\cdot \Psi(2^{-k-1}e_n).
\end{aligned}
\end{equation*}
Without loss of generality, we assume that \textbf{Case 1} holds.

Set $r=2^{-k}$ and define
\begin{equation*}
\omega'=\left\{x\in \omega: |x-e_n|< \mathrm{dist}(e_n,\partial \omega)/2\right\}.
\end{equation*}
By the Harnack inequality and noting that
\begin{equation*}
(u-b_k\Psi)(re_n/2)\geq \frac{a_k-b_k}{2}\cdot \Psi(re_n/2),
\end{equation*}
we have
\begin{equation*}
  u-b_k\Psi\geq C(a_k-b_k)\Psi(re_n/2)~\mbox{ in}~E(\omega',r/4,3r/4).
\end{equation*}

Let $y=x/r$ and
\begin{equation*}
\tilde{u}(y)=\frac{u(x)-(b_k+\mu (a_k-b_k))\Psi(x)}{r^{\alpha_{\omega}}}
\cdot \frac{C}{2(a_k-b_k)\Psi(e_n/2)},
\end{equation*}
where $0<\mu<C/2$ is a small constant to be specified later. Then $\tilde{u}$ satisfies
\begin{equation*}
  \left\{\begin{aligned}
    &M^-(D^2\tilde{u})\leq 0&& ~~\mbox{in}~~E(\omega,1/8,1);\\
    &\tilde u\geq 1&& ~~\mbox{in}~~ E(\omega',1/4,3/4);\\
    &\tilde u= 0&& ~~\mbox{on}~~\partial C_{\omega}\cap (B_1\backslash \bar B_{1/8});\\
    &\tilde u\geq -C\mu/(2\Psi(e_n/2))&& ~~\mbox{on}~~\partial (B_1\backslash \bar B_{1/8})\cap C_{\omega}.
  \end{aligned}\right.
\end{equation*}
Hence, by taking $\mu$ small enough, $\tilde{u}$ satisfies the conditions of \Cref{le.3.1}. Therefore,
\begin{equation*}
 \tilde u\geq 0~~\mbox{ in}~~E(\omega,1/4,3/4).
\end{equation*}
That is,
\begin{equation*}
u\geq \left(b_k+\mu(a_k-b_k)\right)\Psi~~\mbox{ on}~~\partial B_{r/2}\cap C_{\omega}.
\end{equation*}
By the comparison principle,
\begin{equation*}
u\geq \left(b_k+\mu(a_k-b_k)\right)\Psi~~\mbox{ in}~~C_{\omega}^{r/2}.
\end{equation*}

%
%

Let $a_{k+1}=a_k$ and $b_{k+1}=b_k+\mu(a_k-b_k)$. Then
\begin{equation*}\label{e.reg-cone-ak}
a_{k+1}-b_{k+1}=(1-\mu)(a_k-b_k).
\end{equation*}
That is, \cref{e.reg-cone-M} holds for $k+1$. By induction, the proof is completed. \qed~\\

As on the upper half-space, if we have the control of the behavior of $u$ at infinity, the Liouville theorem-type I \Cref{t.reg-co} holds based on the boundary pointwise regularity.

\noindent \textbf{Proof of \Cref{t.reg-co}.} The proof is similar to that of \Cref{co2.1}. By the boundary pointwise regularity \Cref{t.reg-F-cone}, there exists $a\in R$ such that for any $R>0$,
\begin{equation}\label{e3.8}
|u(x)-a\Psi(x)|\leq C\frac{|x|^{\alpha_{\omega}+\alpha}}{R^{\alpha_{\omega}+\alpha}}
\|\Psi\|_{L^{\infty}(\omega)}\|u\|_{L^{\infty}(C^R_{\omega})},~~\forall x\in C^{R/2}_{\omega}.
\end{equation}
Since $u=O(\Psi)$ as $|x|\rightarrow \infty$, for any $R>0$,
\begin{equation*}
\|u\|_{L^{\infty}(C^R_{\omega})}\leq C\|\Psi\|_{L^{\infty}(C^R_{\omega})}\leq CR^{\alpha_{\omega}}.
\end{equation*}
Fix $x\in C_{\omega}$ and let $R\rightarrow \infty$ in \cref{e3.8}. Then $u(x)= a\Psi(x)$. Hence,
\begin{equation*}
  u\equiv a\Psi~\mbox{ in }~C_{\omega}.
\end{equation*}
Obviously, $a=u(e_n)/\Psi(e_n)$. If $u=o(\Psi)$, $a$ must be zero. That is, $u\equiv 0$ in $C_{\omega}$.~\qed~\\

%


In the following, we prove the Liouville theorem-type II for nonnegative solutions. We need the following Carleson type estimate on cones, which can be proved similarly as \Cref{l.cal}.
\begin{lemma}[\textbf{Carleson type estimate on cones}]\label{l.cal-cone}
Let $u(e_n/2)=1$ and $u\geq 0$ be a viscosity solution of
\begin{equation}\label{e.cal-co1}
\left\{\begin{aligned}
&u\in S(\lambda,\Lambda,0)&& ~~\mbox{in}~~C^1_{\omega};\\
&u=0&& ~~\mbox{on}~~(\partial C_{\omega})^{1}.
\end{aligned}\right.
\end{equation}
Then
\begin{equation*}\label{e.cal-co2}
  \|u\|_{L^{\infty}(C^{1/2}_{\omega})}\leq C,
\end{equation*}
where $C$ is universal.
\end{lemma}

Now, we can give the

\noindent \textbf{Proof of \Cref{t.ful-cone}.} By \Cref{t.reg-F-cone}, there exists a constant $a$ such that
for any $R>0$,
\begin{equation}\label{e.t.ful-reg-p}
|u(x)-a\Psi(x)|\leq C\frac{|x|^{\alpha_{\omega}+\alpha}}{R^{\alpha_{\omega}+\alpha}}
\|\Psi\|_{L^{\infty}(\omega)}\|u\|_{L^{\infty}(C^R_{\omega})},~~\forall x\in C^{R/2}_{\omega},
\end{equation}
By the Carleson type estimate on cones,
\begin{equation}\label{e3.9}
  \|u\|_{L^{\infty}(C^R_{\omega})}\leq Cu(Re_n/2).
\end{equation}
From the Hopf lemma \Cref{l.bas-cone},
\begin{equation*}
u(x)\geq C u(Re_n/2)\cdot \Psi(x/|x|)\cdot\frac{|x|^{\alpha_{\omega}}}{R^{\alpha_{\omega}}},
~~\forall~~x\in C^{R/2}_{\omega}.
\end{equation*}
By taking $x=e_n/2$ and we have
\begin{equation}\label{e3.10}
  u(Re_n/2)\leq C u(e_n/2)\cdot R^{\alpha_{\omega}}.
\end{equation}
By combining \cref{e.t.ful-reg-p}, \cref{e3.9} and \cref{e3.10}, we have
\begin{equation*}
|u(x)-a\Psi(x)|\leq C\frac{|x|^{\alpha_{\omega}+\alpha}}{R^{\alpha}}
,~~\forall x\in C^{R/2}_{\omega},
\end{equation*}
Fix any $x\in C_{\omega}$ and let $R\rightarrow \infty$ in above inequality. Then $u(x)=a\Psi(x)$ and $a=u(e_n)/\Psi(e_n)$. Hence,
\begin{equation*}
  u\equiv \frac{u(e_n)}{\Psi(e_n)}\Psi~~~\mbox{ in}~~C_{\omega}.
\end{equation*}
\qed~\\


%
\bibliographystyle{model4-names}
\bibliography{PDE}

\end{document}